\documentclass[12pt]{article}
\usepackage{graphicx}
\usepackage{amssymb}
\usepackage{amsmath}
\numberwithin{equation}{section}

\begin{document}
\textbf{Effective Construction of a Positive Operator which
 does not admit  Triangular Factorization }\\
\begin{center}{\textbf{Sakhnovich Lev}}\end{center}
\emph{99 Cove  ave., Milford, CT, 06461, USA}\\
 E-mail: lsakhnovich@gmail.com\\
 \textbf{Abstract.} We have constructed a concrete example of a
 non-factorable positive operator. As a result, for the
 well-known problems (Ringrose, Kadison and Singer problems) we  replace  existence
 theorems by concrete examples.\\
 \textbf{Mathematical Subject Classification (2000).}  Primary 47A68; Secondary 47A05,
 47A66.\\
 \textbf{Keywords.} Triangular operators, nest algebra, multiplicity 1,  hyperintransitive
 operator.
 \newpage.
\section{Introduction}We introduce the main notions of the
triangular factorization (see  [3, 5, 7] and [10, 11, 16]).\\
In the Hilbert space $L_{m}^{2}(a,b),\quad
(-\infty{\leq}a<b{\leq}\infty)$ we define the orthogonal projectors
\\$P_{\xi}f=f(x) ,\, \, a{\leq}x<\xi \quad $ and
 $\quad P_{\xi}f=0 ,\,\, \xi<x{\leq}b$ , where $f(x)\in L_{m}^{2}(a,b)$.
 
\textbf{Definition 1.1}. A bounded operator $S_{-}$ on
 $L_{m}^{2}(a,b)$  is called lower triangular if for
every
 $\xi$ the relations \begin{equation}
  S_{-}Q_{\xi}=Q_{\xi}S_{-}Q_{\xi} , \end{equation} where
  $Q_{\xi}=I-P_{\xi}$, are true. A bounded operator  $S_{+}$ is called
  upper triangular, if  $S_{+}^*$ is lower triangular.
  
  \textbf{Definition 1.2}.  A bounded, positive definite and invertible
operator $S$ on $L_{m}^{2}(a,b)$ is said to admit the left (the right)
triangular factorization if it can be represented in the form
\begin{equation}
S=S_{-}S_{-}^* \qquad (S=S_{-}^*S_{-}),\end{equation} where $S_{-}$ and $S_{-}^{-1}$ are
lower triangular, bounded
operators.\\
In paper [16] (p. 291)  we formulated the necessary and sufficient conditions
under which a positive definite operator $S$ admits  the triangular
factorization. The factorizing operator $S_{-}^{-1}$ was
constructed in the explicit form.\\
Let us introduce the notations
\begin{equation} S_{\xi}=P_{\xi}SP_{\xi},\quad (f,g)_{\xi}=\int_{0}^{\xi}f^*(x)g(x)dx, \quad f,g{\in}L_{m}^{2}(0,b).\end{equation}

 \textbf{Theorem 1.1} ([16], p.291)  \emph{Let the
bounded and invertible operator S on $L_{m}^{2}(a,b)$ be positive
definite.
For the operator S to admit the left triangular factorization it is
necessary and sufficient that  the following assertions be
true.\\
1.There exists  such a $m{\times}m$ matrix function $F_{0}(x)$ that
\begin{equation}
Tr\int_{a}^{b}F_{0}^*(x)F_{0}(x)dx<\infty,\end{equation} and the
$m{\times}m$ matrix function
\begin{equation} M(\xi)=(F_{0}(x),S_{\xi}^{-1}F_{0}(x))_{\xi}
\end{equation}}
\emph{is absolutely continuous and almost everywhere
\begin{equation} {\mathrm{det}}M'(\xi){\ne}0.\end{equation} 2. The vector functions
\begin{equation}
\int_{a}^{x}v^*(x,t)f(t)dt \end{equation} are absolutely continuous.
Here $f(x){\in}L_{m}^{2}(a,b)$, \begin{equation}
v(\xi,t)=S_{\xi}^{-1}P_{\xi}F_{0}(x) ,
\end{equation}
(In (1.8), the operator $S_{\xi}^{-1}$  is applied columnwise.)\\ 3.
The operator
\begin{equation}Vf=[R^*(x)]^{-1}\frac{d}{dx}
\int_{a}^{x}v^*(x,t)f(t)dt,
\end{equation}
is bounded, invertible and lower triangular  together with its inverse
$V^{-1}$. Here $R(x)$ is a $m{\times}m$ matrix function such that}
\begin{equation}
 R^*(x)R(x)=M'(x).\end{equation}
 \textbf{Corollary 1.1.} ([16], p.293) \emph{If the conditions of Theorem 1.1 are
  fulfilled then the corresponding operator $S^{-1}$ can be represented in the form}
 \begin{equation} S^{-1}=V^*V.\end{equation}
 The formulated theorem allows us to prove that a wide class of
 operators admits the triangular factorization [16].
 
 D.Larson proved  [7] the \emph{existence}  of positive definite and invertible
 but
non-factorable operators.

 In the present
 article we  construct a \emph{concrete example} of positive definite
 and invertible but
 non-factorable operator.
 
 Using this result
  we have managed to substitute  existence
 theorems by concrete examples in the
 well-known problems (Ringrose, Kadison and Singer problems).
 
 \section{Special class of operators}
 We consider the operators $S$ of the form
 \begin{equation}Sf=f(x)-{\mu}\int_{0}^{\infty}h(x-t)f(t)dt,\quad
 f(x){\in}L^{2}(0,\infty),\end{equation}
 where $\mu=\overline{\mu}$ and  $h(x)$ has the representation
 \begin{equation}h(x)=\frac{1}{2\pi}\int_{-\infty}^{\infty}e^{ix\lambda}\rho(\lambda)d\lambda.
 \end{equation}We suppose that the function $\rho(\lambda)$
 satisfies the following conditions\\
 1.
  $\overline{\rho(\lambda)}=\rho(\lambda){\in}L(-\infty,\infty)$.\\
2.The function $\rho(\lambda)$ is bounded, i.e.
 $|\rho(\lambda)|{\leq}M,\quad
  -\infty<\lambda{<}\infty.$\\
  Hence  the function $h(x)\quad
  (-\infty<x<\infty)$ is continuous and the operator
  \begin{equation}Hf=\int_{0}^{\infty}h(x-t)f(t)dt \end{equation}
  is self-adjoint and bounded
  $||H||{\leq}M.$
  We introduce the operators
   \begin{equation}S_{\xi}f=f(x)-{\mu}\int_{0}^{\xi}h(x-t)f(t)dt,\quad
 f(x){\in}L^{2}(0,\xi),\quad 0<\xi<\infty.\end{equation}If
 $|\mu|{\leq}1/M,$
 then the operator $S_{\xi}$ is positive definite ,bounded and
 invertible. Hence we have
 \begin{equation}S_{\xi}^{-1}f=f(x)+\int_{0}^{\xi}R_{\xi}(x,t,\mu)f(t)dt.\end{equation}
 The function $R_{\xi}(x,t,\mu)$ is jointly continuous to the
 variables $x,t,\xi , \mu$.
 M.G.Krein (see [4], Ch.IV , section 7) proved that
 \begin{equation}S_{b}^{-1}=(I+V_{+})(I+V_{-}),\quad
 0<b<\infty,\end{equation}
 where the operators $V_{+}$ and $V_{-}$ are defined in $L^{2}(0,b)$ by the
 relations
\begin{equation}V_{+}^{\star}f=V_{-}f=\int_{0}^{x}R_{x}(x,t,\mu)f(t)dt.
\end{equation}
\textbf{Remark 2.1} The Krein's result (2.6) is true for the Fredholm class of operators.
The operator $S_{b}$ belongs to this class but the operator $S$ does not belong.
When considering the operator $S$ we use Theorem 1.1.
\section{Main example}
Let us consider the operator $S$ ,defined by formula (2.1), where
\begin{equation}h(x)=\frac{\mathrm{sin}{\pi}x}{{\pi}x}=\frac{1}{2\pi}\int_{-\pi}^{\pi}e^{ix\lambda}d{\lambda},\quad
0<{\mu}<1.\end{equation}
In case (3.1) the function $\rho(\lambda)$ has the form
\begin{equation}\rho(\lambda)=1 ,\quad \lambda{\in}[-\pi,\pi];\quad
\rho(\lambda)=0 \quad
 \lambda{\notin}[-\pi,\pi].\end{equation}Hence we have
 $M=1$ and the following statement is true.\\
 \textbf{Proposition 3.1.} \emph{The operator}
 \begin{equation}Sf=f(x)-{\mu}\int_{0}^{\infty}\frac{\mathrm{sin}{\pi}(x-t)}{{\pi}(x-t)}f(t)dt,\quad
 f(x){\in}L^{2}(0,\infty),\quad 0<\mu<1\end{equation}
  \emph{is self-adjoint, bounded, invertible and
 positive definite.}\\
 The following assertion is  the main result of this paper.\\
 \textbf{Theorem 3.1.} \emph{ The bounded  positive definite and invertible operator S, defined by formula (3.3),
does not admit the left triangular factorization}\\
We shall prove Theorem 3.1 by parts.The key results will be written
in the form of lemmas and propositions.\\
Let us  consider the following
 functions
 \begin{equation}Q_{\xi}(\xi,\xi,\mu)=R_{2\xi}(2\xi,2\xi,\mu),\quad
Q_{\xi}(\xi,-\xi,\mu)=R_{2\xi}(2\xi,0,\mu).\end{equation} We use the
relation (see [19], p.16, formula (77))
\begin{equation}
\frac{d}{dt}[Q_{t}(t,t,\mu)]=2Q_{t}^{2}(-t,t,\mu),
\end{equation} and the asymptotic representation (see [9], p.189, formulas (1.16) and (1.17))\begin{equation}
\sigma(x,\mu)=a(\mu)x+b(\mu)+F_{-1}(x)/x + O(1/x^{2}),\quad x{\to}\infty.
\end{equation}where
\begin{equation}a(\mu)=\frac{1}{\pi}\mathrm{log}(1-\mu),\quad
b(\mu)=\frac{1}{2}a^{2}(\mu),\end{equation}
\begin{equation}F_{-1}(x)=\frac{1}{4}a(\mu)^{2}\mathrm{sin}[2x+x_{0}+k\mathrm{log}(x)]+m.\end{equation}
Here the fixed numbers $x_{0},\,k,\, m$ are real and $k{\ne}0.$
The functions $Q_{t}(t,t,\mu)$ and $\sigma(x,\mu)$ are connected by the relation
(see [9], p.189, formula (1.18))
\begin{equation}\sigma(x,\mu)=-2tQ_{t}(t,t,\mu),\quad where \quad
x=2{\pi}t.\end{equation} It follows from (3.6) that
\begin{equation}
\frac{d\sigma(x.\mu)}{dx}=a(\mu)+x^{-1}a(\mu)^{2}\cos[2x+x_{0}+k\mathrm{log}(x)]+O(1/x^{2}),\quad x{\to}\infty.\end{equation}
From (3.6) and (3.9) we deduce that
\begin{equation}Q_{t}(t,t,\mu)=-a(\mu)\pi-b(\mu)/(2t)-F_{-1}(2{\pi}t)/(2{\pi}t^{2})+O(1/t^{3}).\end{equation}
Relations
(3.5), (3.8) and (3.11) imply
\begin{equation}Q^{2}_{t}(-t,t,\mu)=a^{2}(\mu) \mathrm{sin}^{2}[2{\pi}t+x_{0}+k\mathrm{log}(t)]/(4t^{2}) +O(1/t^{3}).\end{equation}
 From (3.4) and (3.12) we deduce the assertion.\\
\textbf{Lemma 3.1} \emph{One of the two relations}
\begin{equation}R_{t}(t,0,\mu)=\epsilon\frac {a(\mu)}{t}\mathrm{sin}[{\pi}t+x_{0}+k\mathrm{log}(t)] +O(1/t^{3/2}),\quad t{\to}\infty
\end{equation}\emph{is valid. Here} $\epsilon=\pm{1}$.\\
 Let us introduce the function
\begin{equation}q_{1}(x)=1+\int_{0}^{x}R_{x}(x,t,\mu)dt.
\end{equation}
According to the
well-known Krein's formula ([4],Ch.IV,formulas (8.3) and (8.14)) we have
\begin{equation}q_{1}(x)=\mathrm{exp}[\int_{0}^{x}R_{t}(t,0,\mu)dt].\end{equation}
\textbf{Lemma 3.2} \emph{The following relation}
\begin{equation}\lim_{x{\to}\infty}{q_{1}(x)}=\frac{1}{\sqrt{1-\mu}} \end{equation}
\emph{is true.}\\
\emph{Proof.} Using the asymptotical formula of confluent hypergeometric function
(see [1], section 6.13) we deduce that the integral
\begin{equation}\mathrm{exp}[\int_{0}^{x}R_{t}(t,0,\mu)dt]=C\end{equation}
converges.
Together with $q_{1}(x)$ we shall consider the function
\begin{equation}q_{2}(x)=M(x)+\int_{0}^{x}M(t)R_{x}(x,t,\mu)dt,
\end{equation}where
\begin{equation}M(x)=\frac{1}{2}-{\mu}\int_{0}^{x}\frac{sin(s\pi)}{s\pi}ds.\end{equation}
The function $M(x)$ can be represented in the form
\begin{equation}M(x) =(1-\mu)/2+q(x).\end{equation}
Using asymptotic of sinus integral (see [2],Ch.9,formulas (2) and (10)), we have
\begin{equation}q(x)=O(1/x),\quad x{\to}\infty.\end{equation}
From (3.15), (3.17) and (3.20) we
deduce
\begin{equation}\lim_{x{\to}\infty}{q_{1}(x)}=C,\,\lim_{x{\to}\infty}{q_{2}(x)}=C(1-\mu)/2.
\end{equation}
Taking into account the relation
 $q_{1}(x)q_{2}(x)=1/2$ (see [15], formulas (53), (56))  we deduce the equality
\begin{equation}\lim_{x{\to}\infty}{q_{2}(x)}=1/(2C)\end{equation}
Relation (3.16) follows directly from formulas (3.22),(3.23) and inequality $C>0$.
The lemma is proved.\\
We note, that relations (3.15) and (3.16) define the sign of $\epsilon$  in equality (3.13).\\
The functions $q_{1}(x)$ and $q_{2}(x)$ generate the $2{\times}2$ differential system
\begin{equation}\frac{dW}{dx}=izJH(x)W,\quad W(0,z)=I_{2}.\end{equation}Here $W(x,z)$,  $J$, $H(x)$
are   $2{\times}2$ matrix functions and
\begin{equation} H(x)=\left[\begin{array}{cc}
                       q_{1}^{2}(x) & 1/2 \\
                       1/2 & q_{2}^{2}(x)
                     \end{array}\right],
                     J=\left[\begin{array}{cc}
                        0 & 1 \\
                        1 & 0
                      \end{array}\right].\end{equation}
It is easy to see that
\begin{equation}JH(x)=T(x)PT^{-1}(x),\end{equation}
where
\begin{equation}T(x)=\left[\begin{array}{cc}
                             q_{2}(x) & -q_{2}(x) \\
                             q_{1}(x) & q_{1}(x)
                           \end{array}\right],
P=\left[\begin{array}{cc}
    1 & 0 \\
    0 & 0
  \end{array}\right].\end{equation}
Let us consider the matrix function
\begin{equation}V(x,z)=e^{-ixz/2}T^{-1}(x)W(x,z)T(0).\end{equation}
Due to (3.24)-(3.28) we have
\begin{equation}\frac{dV}{dx}=(iz/2)jV-Q(x)V, \quad V(0)=I_{2},\end{equation}
where
\begin{equation}Q(x)=\left[\begin{array}{cc}
                             0 & B(x) \\
                             B(x)& 0
                           \end{array}\right],
j=\left[\begin{array}{cc}
          1 & 0 \\
          0 & -1
        \end{array}\right],
\end{equation}
\begin{equation} B(x)=\frac{q_{1}^{\prime}(x)}{q_{1}(x)}=R_{x}(x,0,\mu).\end{equation}
Let us introduce the functions
\begin{equation}\Phi_{n}(x,z)=v_{1,n}(x,z)+v_{2,n}(x,z),\quad (n=1,2),\end{equation}
\begin{equation}\Psi_{n}(x,z)=i[v_{1,n}(x,z)-v_{2,n}(x,z)], \quad (n=1,2),\end{equation}
where $v_{i,n}(x,z)$ are elements of the matrix function $V(x,z)$.
It follows from (3.29) that
\begin{equation}\frac{d\Phi_{n}}{dx}=(z/2)\Psi_{n}-B(x)\Phi_{n},\quad \Phi_{1}(0,z)=\Phi_{2}(0,z)=1 \end{equation}
\begin{equation}\frac{d\Psi_{n}}{dx}=-(z/2)\Phi_{n}+B(x)\Psi_{n}
\quad \Psi_{1}(0,z)=-\Psi_{2}(0,z)=i .\end{equation}
Let us consider again the differential system (3.24) and the  solution $W(x,z)$ of this system.
The element $w_{2,1}(\xi,z)$ of the matrix function $W(x,z)$
can be represented in the form  (see [13], p.54, formula (2.5))
\begin{equation}w_{2.1}(\xi,z)=iz((I-Az)^{-1}1,S_{\xi}^{-1}1)_{\xi},\end{equation}
where the operator $A$ has the form
\begin{equation}Af=i\int_{0}^{x}f(t)dt.\end{equation}
It is well-known that
\begin{equation}(I-Az)^{-1}1=e^{izx}.\end{equation}
Now we need the relations (see [12], Ch.1, formulas (1.37) and (1.44)):
\begin{equation}S_{\xi}1=M(x)+M(\xi-x),\quad S_{\xi}=U_{\xi}S_{\xi}U_{\xi}, \end{equation}
where the function $M(x)$ is defined by relation (3.20) and\\ $U_{\xi}f(x)=\overline{f(\xi-x)},\quad 0{\leq}x{\leq}\xi.$
  It follows from (3.39) that
\begin{equation}S_{\xi}1=1-\mu +q(x)+U_{\xi}q(x).\end{equation}Hence the relation
\begin{equation}S_{\xi}^{-1}1=\frac{1}{(1-\mu)}[1-R_{\xi}(x)-U_{\xi}R_{\xi}(x)],\quad q(x)=O(1/x),\end{equation}
where $S_{\xi}^{-1}q(x)=R_{\xi}(x).$
is true. According to formulas (3.36) and (3.41) the following statement is true.\\
\textbf{Lemma 3.3} \emph{The function $w_{2.1}(\xi,z)$ has the form}
\begin{equation}w_{2.1}(\xi,z)=e^{iz\xi}G(\xi,z)-\overline{G(\xi,\overline{z})},\end{equation}
\emph{where}
\begin{equation} G(\xi,z)=\frac{1}{1-\mu}[1-iz\int_{0}^{\xi}e^{-izx}R_{\xi}(x)dx].\end{equation}
We can obtain another representation of $w_{2.1}(\xi,z)$ without using the operator $S_{\xi}^{-1}$.
Indeed, it follows from
(3.28) and (3.32),(3.33) that
\begin{equation}W(x,z)=(1/2)e^{ixz/2}T(x)\left[\begin{array}{cc}
                                            \Phi_{1}-i\Psi_{1}& \Phi_{2}-i\Psi_{2} \\
                                         \Phi_{1}+i\Psi_{1} & \Phi_{2}+i\Psi_{2}
                                          \end{array}\right]T^{-1}(0).\end{equation}
According to equality (3.14) we have $q_{1}(0)=1$. Due to (3.27) we infer
\begin{equation}T(0)=\left[\begin{array}{cc}
                             1/2& -1/2 \\
                             1 & 1
                           \end{array}\right],\quad T^{-1}(0)=\left[\begin{array}{cc}
                             1& 1/2 \\
                             -1 & 1/2
                           \end{array}\right].\end{equation}
It follows from (3.16),(3.22) and (3.27) that
 \begin{equation}T(x){\to}\left[\begin{array}{cc}
                            1/(2C) & -1/(2C) \\
                            C & C
                          \end{array}\right],\quad x{\to}\infty,\quad C=1/\sqrt{(1-\mu)}.
                          \end{equation}
Hence in view of (3.45)-(3.46) the following assertion is true.\\
\textbf{Lemma 3.4.} \emph{ The function $w_{2.1}(x,z)$ has the form ($x{\to}\infty$ )}
\begin{equation}w_{2.1}(x,z)=Ce^{ixz/2}\phi(x,z)(1+o(1)),\,\phi(x,z)=\Phi_{1}(x,z)-\Phi_{2}(x,z).
\end{equation}\\
Further we plan to use one Krein's result [6]. To do it we introduce the functions
\begin{equation}P(x,z)=e^{ixz/2}[\Phi(x,z)-i\Psi(x,z)]/2,\end{equation}
\begin{equation}P_{\star}(x,z)=e^{ixz/2}[\Phi(x,z)+i\Psi(x,z)]/2,\end{equation}
where
\begin{equation}\Phi(x,z)=\Phi_{1}(x,z)+\Phi_{2}(x,z),\,\Psi(x,z)=\Psi_{1}(x,z)+\Psi_{2}(x,z).
\end{equation}
Using (3.34),(3.35) and (3.48),(3.40) we see that $P(x,z)$ and $P_{\star}(x,z)$ is the solution
of the following Krein's system
\begin{equation}\frac{dP}{dx}=(iz/2)P-B(x)P_{\star},\quad  \frac{dP_{\star}}{dx}=-B(x)P,
\end{equation} where
\begin{equation}P(0,z)=P_{\star}(0,z)=1
\end{equation}The coefficient $B(x)$ belongs to $L^{2}(0,\infty)$ (see (3.13) and (3.31)).
Hence  the following Krein's results are true [6].\\
\textbf{Proposition 3.2} 1)\emph{ There exists the limit}
\begin{equation}\Pi(z)=\lim_{x{\to}\infty} P_{\star}(x,z),\end{equation}
\emph{where the convergence is uniform at any bounded closed set $z$ of the open half-plane
$Im z>0$. }\\
2)\emph{The function $\Pi(z)$ can be represented in the form}
\begin{equation} \Pi(z)=\frac{1}{\sqrt{2\pi}}\exp[\frac{1}{2i\pi}\int_{-\infty}^{\infty}
\frac{1+tz/2}{(t-z/2)(1+t^{2})}(\log{\sigma}^{\prime}(t))dt+i\alpha],\end{equation}
\emph{where} $\alpha=\overline{\alpha}$.\\ Here  $\lambda=z/2$ is the spectral parameter
of system (3.52),
 $\sigma(u)$ is the spectral function of this system and
  is defined by the relation (see [6],formula (2)):
\begin{equation}\sigma^{\prime}(u)=1/(2\pi)\quad (|u|>\pi),\quad
 \sigma^{\prime}(u)=(1-\mu)/(2\pi)\quad (|u|<\pi).\end{equation}
In case (3.3) the following conditions are fulfilled:
\begin{equation} 1-\mu{\leq}||S||{\leq}1+\mu,\quad \int_{0}^{\infty}|h(x)|^{2}dx<\infty.\end{equation}
Therefore in case (3.3) we can use Proposition 1 and formula (2.15) from paper [17].
Hence in formula (3.54) we have
\begin{equation} \alpha =0.\end{equation}
It follows from (3.55) and (3.57) that  (3.54) takes the form
\begin{equation} \Pi(z)=\exp{[\frac{1}{2i\pi}\int_{-\pi}^{\pi}\frac{\log(1-\mu)}{t-z/2}dt]}\end{equation}
Now we need the relation
\begin{equation}\int_{-\pi}^{\pi}\frac{1}{t-z/2}dt=\frac{1}{2}\log\frac{(2\pi-x)^{2}+y^{2}}{(2\pi+x)^{2}+y^{2}}+iL(x,y),\end{equation}
where $z=x+iy$ and
\begin{equation}L(x,y)=\arctan\frac{2\pi-x}{y}+\arctan\frac{2\pi+x}{y}.\end{equation}
Due to (3.50) and (3.60) we obtain that
\begin{equation}\lim_{y{\to}+0}\int_{-\pi}^{\pi}\frac{1}{t-z/2}dt=
\log|\frac{2\pi-x}{2\pi+x}|
+i{\pi}{\chi(x)},\end{equation}where $\chi(x)=1$ when $|x|<2\pi$, $\chi(x)=0$ when $|x|>2\pi.$
Taking into account equalities (3.58) and  (3.61)  we have
\begin{equation}\Pi(x)=\lim_{y{\to}+0}\Pi(x+iy)=\left|\frac{2\pi-x}{2\pi+x}\right|^{1/(2i\pi)}\sqrt{1-\mu\chi(x)},\end{equation}
\begin{equation}\Pi(0)=\lim_{y{\to}+0}\Pi(iy)=\sqrt{1-\mu}.\end{equation}
Using the spectral function $\sigma(u)$ we can construct the  Weyl-Titchmarsh function
$v(z)$
of system (3.51) (see [13]):
\begin{equation} v(z)=\int_{-\infty}^{\infty}
\frac{1+tz/2}{(t-z/2)(1+t^{2})}{\sigma}^{\prime}(t)dt+\frac{2z\mu}{4{\pi}^{2}-z^{2}}.\end{equation}
It follows from  (3.59) that
\begin{equation}\lim_{y{\to}+0}v(z)=i\frac{1-\mu}{2}-\frac{\mu}{\pi}\log\frac{|2\pi-x|}{|2\pi+x|}+\frac{2x\mu}{4{\pi}^{2}-x^{2}}
,\quad x{\in}(-2\pi.2\pi).\end{equation} Hence the relation
\begin{equation}v(0)=\lim_{y{\to}+0}v(iy)=i(1-\mu)/2\end{equation}is true. Together with
$P(x,z)$ and $P_{\star}(x,z)$ let us introduce another solution $\hat{P}(x,z)$ and $\hat{P}_{\star}(x,z)$
of system (3.51),where
\begin{equation}\hat{P}(0,z)=1/2,\quad \hat{P}_{\star}(0,z)=-1/2.\end{equation}The following statements are proved in the book
([13],Ch.10).\\
\textbf{Proposition 3.3} 1)\emph{ There exist the limits}
\begin{equation}\hat{\Pi}(z)=\lim_{x{\to}\infty} \hat{P}_{\star}(x,z),\quad \lim_{x{\to}\infty} \hat{P}(x,z)=0.\end{equation}
\emph{where the convergence is uniform at any bounded closed set $z$ of the open half-plane
$Im z>0$. }\\
2)\emph{The equality holds}
\begin{equation}v(z)=-i{\Pi}^{-1}(z)\hat{\Pi}(z),\quad \Im{z}>0.\end{equation}
According to (3.60),(3.63) and (3.66) the relations
\begin{equation}\hat{\Pi}(z)=iv(z)\Pi(z),\quad \hat{\Pi}(0)=-[(1-\mu)^{3/2}]/2.\end{equation}
are valid.
We introduce the function
\begin{equation}\psi(x,z)=\Psi_{1}(x,z)-\Psi_{2}(x,z).\end{equation}The functions $\phi(x,z)-i\psi(x,z)$ and
$\phi(x,z)+i\psi(x,z)$ satisfy the the equation (3.51) and $\phi(0,z)=0,\, \psi(0.z)=2i.$ Hence we have
\begin{equation}\hat{P}(x,z)=\frac{1}{4}[\phi(x,z)-i\psi(x,z)]e^{ixz/2},\quad
\hat{P}_{\star}(x,z)=\frac{1}{4}[\phi(x,z)+i\psi(x,z)]e^{ixz/2}.\end{equation}
It follows from (3.72) that
\begin{equation}\hat{P}(x,z)+\hat{P}_{\star}(x,z)=\frac{1}{2}\phi(x,z)e^{ixz/2}.\end{equation}
\emph{Proof of Theorem 3.1.} Let us suppose that the operator $S$ admits the factorization.
According to Theorem 1.1 the corresponding operator $S_{-}^{-1}$ has the form
$S_{-}^{-1}=I+V_{-}$, where $V_{-}$ is defined by (2.7). Hence the the operator function
 $S_{\xi}^{-1}$ strongly converges to the operator  $S^{-1}$ when $\xi{\to}\infty$.Then
 the function $R_{\xi}(x)=S_{\xi}^{-1}q(x)$ strongly converges to $R(x)=S^{-1}q(x)$ and
 $R(x){\in}L^{2}(0,\infty)$. Using (3.42) and (3.43) we have
 \begin{equation} e^{-iz\xi/2}w_{2.1}(\xi,z)=e^{iz\xi/2}G(z)-e^{-iz\xi/2}\overline{G(z)}+o(1),\,
 \xi{\to}\infty,\end{equation}
where
\begin{equation} G(z)=\frac{1}{1-\mu}[1-iz\int_{0}^{\infty}e^{-izx}R(x)dx], \,z=\overline{z}.
\end{equation}
Taking into account  Lemma 3.4, Proposition (3.3) and relation (3.73) we have
 \begin{equation} e^{-iz\xi/2}w_{2.1}(\xi,z)=e^{iz\xi/2}H(z)-e^{-iz\xi/2}\overline{H(z)}+o(1),\,
 \xi{\to}\infty,\end{equation}
where
\begin{equation}\overline{H(z)}=-2C\Pi(z), \,z=\overline{z}.\end{equation}
 Comparing formulas (3.75) and (3.77) we see that
\begin{equation}G(0)=1/(1-\mu){\ne}H(0)=(1-\mu).\end{equation}Hence the relation
 (3.75) is not true, i.e. the operator $S$ does not admit the factorization.
 The theorem is proved.
\section{Examples instead of existence theorems} Let the nest $N$ be the family of subspaces $Q_{\xi}L^{2}(0,\infty).$
The corresponding \emph{nest algebra} $Alg{N}$ is the algebra of all linear bounded operators
 in the space $L^{2}(0,\infty)$ for which every member of $N$ is invariant subspace. Let us denote by $D_{N}=Alg(N){\bigcap}Alg(N)^{\star}$.
 The set $N$ has \emph{ multiplicity one} if diagonal $D_{N}$ is abelian, i.e.
 $D_{N}$ is commutative algebra.
 We can see that the lower triangular operators $S_{-}$ form the algebra
 $alg(N)$  , the corresponding diagonal $D_{N}$ is abelian and consists
 of the commutative  operators
 \begin{equation}T_{\phi}f=\phi(x)f,\quad f{\in}L^{2}(0,\infty),\end{equation}
 where $\phi(x)$ is bounded. Hence the introduced nest $N$ has the multiplicity 1.
  We obtain the following concrete answer to Ringrose's question (see [3],[7]).\\
\textbf{Proposition 4.1.} \emph{Let the positive definite ,
invertible operator $S$ is defined by the relation (3.4).
 The  set $S^{1/2}N$ fails to have  multiplicity 1}\\
\emph{Proof.}  We use the well-known result (see [3] ,p.169).\\
\emph{The following assertions are equivalent:}\\
\emph{1. The  positive definite , invertible operator T admits factorization.\\
2. $T^{1/2}$ preserves  the multiplicity of $N$.}\\
(We stress that in our case the set $N=Q_{\xi}L^{2}(0,\infty)$ is fixed.)
The operator $S$  does not admit the factorization.Therefore
the  set $S^{1/2}N$ fails to have  multiplicity 1. The proposition is proved.\\
Let as consider the operator
\begin{equation} Vf=\int_{0}^{x}e^{-(x+y)}f(y)dy, \quad f(x){\in}L^{2}(0,\infty)\end{equation}
An operator is said to be \emph{hyperintransitive} if its lattice of invariant
subspaces contains a multiplicity one nest. We note that the lattice of invariant
subspaces of the operator  $V$ coincides with $N$ (see [8] and [18],  Ch.11, Theorem 150). Hence we  deduce the answer
to Kadison-Singer [5] and to Gohberg and Krein [4] question.\\
\textbf{Corollary 4.1.} \emph{The operator $W=S^{1/2}VS^{-1/2}$  is non-hyperintransitive
compact operator.}\\Indeed the lattice of the invariant subspaces of
 the operator $W$ coincides with $S^{1/2}N.$\\
\textbf{Corollary 4.2.}\emph{ If $\epsilon>0$ and $0<\mu<\epsilon$
 then $S=I+K$, where $||K||<\epsilon$.  The corresponding operator $S$ does not admit
 the factorization.}\\
\textbf{Remark 4.1.} The existence parts of Theorem 3.1, Proposition 4.1 , Corollaries  4.1
and 4.2
 are proved by
D.R.Larson [7].\\
But we have constructed the concrete
 simple examples of the operator $S$ of the convolution type and the operator $V$.\\
 \textbf{Remark 4.2.}The factorization problems
 of the convolution type operators were discussed  in our papers [10],[11] and [16].
\begin{center}{References }\end{center}
1.\textbf{Bateman H. and Erdelyi A.,} \emph{Higher Transcendental Functions}, New York, v.1,
1953.\\
2.\textbf{Bateman H. and Erdelyi A.,} \emph{Higher Transcendental Functions}, New York, v.2.,
1953.\\
3. \textbf{Davidson K.R.,} \emph{ Nest Algebras. Triangular forms for operator algebras on 
Hilbert space.},  Pitman Research Notes in Mathematics Series 191,  New York,  Wiley
(1988).\\
4.\textbf{ Gohberg I. and Krein M.G.,} \emph{Theory and Applications of Volterra
Operators in Hilbert Space,} Amer. Math. Soc.,
Providence, (1970).\\
5. \textbf{Kadison R. and Singer I.,} \emph{Triangular Operator Algebras,}
\emph{Amer. J. Math. 82, pp. 227-259, (1960)}\\
6.\textbf{Krein M.G.,} \emph{Continuous Analogues of Proposition on
Polynomial Orthogonal on the
Unit Circle,} Dokl.Akad.Nauk SSSR, 105 (1955), 637-640 (Russian).\\
7.\textbf{ Larson D.R.,} \emph{Nest Algebras and Similarity Transformation,}
Ann. Math., 125, p. 409-427, (1985).\\
8. \textbf{Livshits M.S.,} \emph{Operators, Oscillations, Waves (Open Systems)},
Amer. Math. Soc., Providence, 1973.\\
9.\textbf{McCoy B. and Tang S.,}\emph{Connection Formulae for Painleve V Functions,} Physica
20D (1986) 187-216.\\
10.\textbf{ Sakhnovich L.A.,} \emph{ Factorization  of Operators in $L^{2}(a,b)$,}
Functional Anal. and Appl., 13, p.187-192,
(1979), (Russian).\\
11.\textbf{ Sakhnovich L.A.,} \emph{Factorization of Operators in $L^{2}(a,b)$,} in
Linear and Complex Analysis, Problem Book, (Havin V.P., Hruscev S.V.
and Nikol'skii N.K. (ed.)) Springer Verlag, 172-174
(1984).\\
12.\textbf{ Sakhnovich L.A.,} \emph{Integral Equations with Difference Kernels
on Finite Intervals,} Operator Theory, Advances and Appl. v.84, 1996.\\
13.\textbf{ Sakhnovich L.A.,} \emph{Spectral Theory of Canonical Differential Systems.
 Method of Operator Identities,} Operator Theory, Advances and Appl. v.107, 1999.\\
14.\textbf{ Sakhnovich L.A.,} \emph{Spectral Theory of a Class of Canonical Systems.}
Func.Anal.Appl.34, 119-128, 2000.\\
15.\textbf{ Sakhnovich L.A.,} \emph{On Reducing the Canonical System to Two
Dual Differential Systems,} Journal of Math. Analysis and Appl., 255, 499-509,
2001.\\
16.\textbf{ Sakhnovich L.A.,} \emph{On Triangular  Factorization  of positive  Operators,}
Operator Theory: Advances
and Appl., vol.179 (2007) 289-308.\\
17.\textbf{ Sakhnovich L.A.,} \emph{On Krein's Differential System and its Generalization,}
Integral Equations and Operator Theory, 55, 561-572, 2006.\\
18. \textbf{Tithchmarsh E.C.,} \emph{Introduction to the Theory of Fourier Integrals,}
Oxford, 1937.\\
19. \textbf{Tracy C.A. and Widom H.,} \emph{Introduction to Random
Matrices,} Springer Lecture Notes in Physics, 424, (1993),103-130 .\\

\end{document}